\documentclass[12pt,reqno]{amsart}

\usepackage{amsfonts}
\usepackage{amssymb}
\usepackage{amsthm}
\usepackage{amsmath}

\newtheorem{theo}{Theorem}[section]
\newtheorem{propo}[theo]{Proposition}

\theoremstyle{definition}

\theoremstyle{remark}

\begin{document}

\title[On the enumeration of some D-optimal designs]
{On the enumeration of some D-optimal designs}
\author[W. P. Orrick]{William P. Orrick}
\address{Department of Mathematics, Indiana University,
Bloomington IN 47405, USA}
\date{\today}

\begin{abstract}
Two matrices with elements taken from the set $\{-1,1\}$ are Hadamard
equivalent if one can be converted into the other by a sequence of
permutations of rows and columns, and negations of rows and
columns.  In this paper we summarize what is known about the 
number of equivalence classes of matrices having maximal determinant.
We establish that there are 7 equivalence classes for matrices of
order 21 and that there are at least 9,884 equivalence classes for
matrices of order 26.  The latter result is obtained primarily using 
a switching technique for producing new designs from old.
\end{abstract}

\maketitle

\section{Introduction}\label{sect:intro}
In this paper we consider exact D-optimal first order saturated designs with $n$ observations, which can be thought of simply as $n\times n$ matrices
with entries $\pm1$ having maximal determinant.   We are concerned with the enumeration
and classification of such designs up to equivalence, which is defined below.  

We have two
main goals.  First is to provide a summary of previously known enumeration results and
proofs of some new results reported
in ref.~\cite{KhaOrr06}.  (This reference can also be consulted for a concise summary of what
is known about maximal determinants.)  To this end, we compile from the literature a list of all known results.  The result for $n=11$ seemed to require
confirmation, which we undertook as part of our task.  Then we examine two orders
that have not yet received any attention, $n=21$ and $n=26$.  In order 21, we
find that there are 7 equivalence classes of designs.  We do not attempt a complete
enumeration in order 26, but we do show that there are at least 9,884 inequivalent designs by
performing complete enumerations of subclasses known as Q-classes.
The large number of designs we find in order 26 is surprising to us considering that there are only 78 designs in order 25, and
487  in order 28. 

Our second goal is to demonstrate the efficacy of a technique known as {\em switching.}  In the
special case of Hadamard matrices, large numbers of designs have been shown to exist by
{\em doubling,}~\cite{LinWalLie93,LamLamTon00,LamLamTon01} which means constructing new Hadamard matrices from Hadamard matrices of
half the size.  Recently, it has been shown that switching has the potential to produce far more
Hadamard matrices than can be produced by doubling, and that switching is also applicable to
orders  which are not multiples of 8 (so that doubling cannot be used)~\cite{Orr05b}.  In this paper, we will
demonstrate that switching can also often be applied to D-optimal designs more generally, at
least when the order is congruent to 0, 1, 2, 3, or $4 \bmod 8$.  Not only does switching produce
far more designs than doubling does, it is also more widely applicable.  Apart from the case of
Hadamard matrices whose size is a multiple of 8, doubling can only be used  when the order
is twice that of a symmetric balanced incomplete block design with parameters $(2q^2+2q+1,q^2,q(q-1)/2)$, that is, when the order is of the form $4q(q+1)+2$.  Order 26 corresponds to $q=2$.  The next
such case is order 50 ($q=3$).  In order 26  we find that only a small fraction of the designs
produced by switching come from a doubling construction.  In the last section, we describe the
application of switching to a wide variety of D-optimal designs.

\section{Preliminaries}\label{sect:prelim}
Two matrices with elements taken from the set $\{-1,1\}$ are {\em Hadamard
equivalent} (or simply {\em equivalent})  if one can be converted into the other by a sequence of
permutations of rows and columns, and negations of rows and
columns.   If $R$ denotes a $\pm1$ matrix of order $n$ and we wish to discover conditions under which $R$ has
maximal determinant,  it is helpful
to consider the Gram matrices  $M=RR^{\mathrm{T}}$ and $M'=R^{\mathrm{T}}R$.

The notion of equivalence for matrices $R$ carries over to equivalence of the corresponding $M$
and $M'$ which are necessarily symmetric.  The appropriate notion of equivalence for the latter involves {\em simultaneous}
permutation of rows and the corresponding columns, and {\em simultaneous} negation of rows
and the corresponding columns.  When $n$ is odd, we may, by normalizing $R$ appropriately,
take all entries in $M$ and $M'$ to be congruent to $n$ mod 4.

If $M=M'=nI$ then $R$ certainly has maximal determinant---it is a Hadamard matrix---but this condition cannot hold unless $n=0$, $1$,
or $n\equiv0\bmod4$.  

The largest possible determinant for $n\equiv1\bmod 4$ occurs when $R$ is
equivalent to a matrix with Gram matrix $M=M'=(n-1)I_n+J_n$~\cite{Bar33}.  Such an $R$ is necessarily
the $\pm1$ incidence matrix of a symmetric balanced incomplete block design with parameters $(k^2+(k+1)^2,k^2,k(k-1)/2)$, which implies that $n=k^2+(k+1)^2$ for some integer $k$.
Such designs are known for  $n=1$, 5, 13, 25, and many higher orders.  For $n=9$, $17$, and $21$ which do not meet the necessary condition
for this best possible form, the optimal Gram matrices have been found by computer search.
In all three cases, $M=M'$ and
they differ from the best possible form only in a single row and column, which we may take to be
the first.  For order 17, the off-diagonal entries in the first row and column all equal $-3$ rather
than 1~\cite{MoyKou82}.  For orders $9$ and $21$, the first row and column contain one 5 and four 5s respectively,
with the remaining off-diagonal entries equal to 1~\cite{EhlZel62,ChaKouMoy87}.

If $n\equiv2\bmod4$ and then the largest possible determinant occurs
when $R$ is equivalent to a matrix with Gram matrices $M=M'=(n-2)I_n+2I_2\otimes J_{n/2}$ where
$J_k$ is the $k\times k$ all 1 matrix~\cite{Ehl64a,Woj64}.  
Writing
\begin{equation*}
R=\begin{bmatrix}A & B\\C & D\end{bmatrix}
\end{equation*}
with submatrices $A$, $B$, $C$, and $D$ of order $n/2$, we may, by permuting
and negating  rows and columns if necessary, take the row
and column sums of $A$ and $-D$ to be  $a$ and the row and column sums of $B$ and $C$ to be $b$
with $a\ge b\ge0$ and $a^2+b^2=2n-2$.  Hence a necessary condition on the existence of a design
with this best possible $M$ is that $2n-2$ be expressible as the sum of two squares~\cite{Ehl64a}.  For $n\le50$
this excludes $n=22$ and 34, but designs do exist for the other orders~\cite{Ehl64a,Yan66}.  The optimal $M$ is not known for $n=22$ or $34$.

For $n\equiv3\bmod4$ the best possible $M$ is more complicated~\cite{Ehl64b}, and is not known to be
achieved in any order.  For $n\le15$, the optimal $M$ has been found by computer
search~\cite{GalKie80,Orr05a}.

\section{Known classification and enumeration results}\label{sect:known}
In 1946 Williamson~\cite{Wil46} found the maximal determinants up to $n=7$ and showed that there is a unique optimal matrix up to equivalence in each order.
In his 1988 dissertation~\cite{Smi88}, Smith states that this can be shown up to $n=10$.  Certainly uniqueness of the optimal design
for each of these orders, and for $n=12$ and $13$, can be established by hand without difficulty, especially when the
proved optimal forms of $M=M'$ are assumed.  For other values of $n$, we provide the following survey of known results.

The best studied case is that of Hadamard matrices, and the numbers of inequivalent matrices in orders 16, 20, 24, and 28 are known to be 5, 3, 60, and 487~\cite{Hal61,Hal65,ItoLeoLon81,Kim89,Kim94b}.  

In the case $n\equiv1\bmod 4$, we have that the number of inequivalent designs in order 25 is 78~\cite{Den82}.  
This result was established by Denniston, who used a switching method (which, as we shall demonstrate, is more generally applicable, in particular to order 26) to generate his list of designs.  He established completeness of this list by another
method.   We will prove that there are 3 designs in order 17 and 7 designs in order 21 up to
equivalence.  The result in order 17 is a simple consequence of results already in the literature.

The case $n\equiv2\bmod4$ has been studied by Cohn~\cite{Coh94} who showed that there is a unique design in order 14, and that there are 3 distinct designs in order 18.  The next size currently amenable to study is order 26 which we investigate here.  We find a lower bound of 9,884 inequivalent matrices.  In order 22, for which even the value of the maximal determinant has not been
established, 30 distinct matrices have been found by means of a gradient ascent algorithm which all
have determinant $195312500\times2^{21}$.  The latter is the current determinant record, which was discovered by Bruce Solomon and collaborators~\cite{OrrSolDowSmi03}.

The case $n\equiv3\bmod4$ is the most difficult.  We discuss $n=11$ below.  There is a unique optimal design in order 15~\cite{Orr05a}.  
The maximal determinant value is not known with certainty in any higher order.  In order $19$ the current record, $3411968\times2^{18}$, is achieved by matrices which have two distinct forms of $M=M'$. One of which was found by Smith~\cite{Smi88}, the other by Cohn.  Orrick and Solomon~\cite{OrrSol05} constructed a third matrix which has the same $M=M'$
as Cohn's matrix.  No other matrices with determinant equalling the record value are known.  In order 23, 14 distinct matrices achieving the current record determinant, $662671875 \times2^{22}$~\cite{OrrSolDowSmi03},
have been obtained by gradient ascent.  In order 27, a new determinant record, $198087192576\times2^{26}$, was recently
established by Hiroki Tamura~\cite{Tam05}.  By switching, 66  equivalence classes have been obtained
using Tamura's design as a seed matrix.

Current knowledge is summarized in the table.  For those values of $n$ for which the maximal determinant
is not known, the current record is listed, with value written in italics.

\begin{table}
\caption{Number of inequivalent matrices achieving the maximal determinant (current record
in italics where maximal determinant is not known).}
\begin{tabular}{ll|ll|ll|ll}
$n$ & $N$ & $n$ & $N$ & $n$ & $N$ & $n$ & $N$\\
\hline
1 & 1 & 2 & 1 & 3 & 1 & 4 & 1\\
5 & 1 & 6 & 1 & 7 & 1 & 8 & 1\\
9 & 1 & 10 & 1 & 11 & 3 & 12 & 1\\
13 & 1 & 14 & 1 & 15 & 1 & 16 & 5\\
17 & 3 & 18 & 3 & {\em 19} & $\mathit{\ge3}$ & 20 & 3\\
21 & 7 & {\em 22} & $\mathit{\ge30}$ & {\em 23} & $\mathit{\ge14}$ & 24 & 60\\
25 & 78 & 26 & $\ge9884$ & {\em 27} & $\mathit{\ge66}$ & 28 & 487\\
\end{tabular}
\end{table}

\section{D-optimal designs of order 11}
It seems probable that Ehlich proved there are 3 inequivalent matrices of maximal determinant.
He certainly
showed that there are 3 inequivalent forms that $M=M'$ can take.   His proof was never published
although an account of it appears in a paper of Galil and Kiefer~\cite{GalKie80}.  Nevertheless, there appears to be no clear
statement anywhere in the literature that each of these 3 forms decomposes as $M=RR^{\mathrm{T}}$
in an essentially unique way.  We have taken the trouble to confirm that this is the case.  

The proof could be done by hand, but would be tedious.  Instead we adapted some existing
programs for the purpose.  The backtracking decomposition program used in~\cite{Orr05a}
to study the cases $n=29$, $33$, and $37$ was run on each of Ehlich's 3 matrices $M$.  Since
our backtracking program does not perform complete isomorph rejection, a set of different matrices
$R$ was produced in each case.  We then used Brendan McKay's {\em nauty}~\cite{McK04,McK79} to establish that
all matrices in each set are equivalent to each other.  Hence there are 3 equivalence classes of
maximal determinant matrices in  order 11.

\section{D-optimal designs of order 17}
In order 17, it was shown by Moyssiadis and Kounias~\cite{MoyKou82} that the maximal determinant matrix $R$ must be equivalent to a matrix whose Gram matrices $M=M'$  have all 16 off-diagonal entries in the first row and in the first column equal to $-3$, while all remaining off-diagonal entries equal 1.  They noted that, up to overall sign, this implies that the first row and column of $R$ consist
entirely of 1s while the $16\times16$ ``core''  is a regular Hadamard matrix of order 16.  The latter
is the $\pm1$ incidence matrix of a $2$-$(16,6,2)$ design.
(A regular Hadamard matrix is one with constant row and columns sums, and must have perfect
square order.)  Husain~\cite{Hus45} proved that there are 3 such designs up to isomorphism.
From this it follows that there are 3 inequivalent D-optimal designs of order 17.


\section{D-optimal designs of order 21}
Chadjipantelis, Kounias, and Moyssiadis found that the optimal matrix $M=M'$ is one with four 5s
in the first row and column, and 1s everywhere else~\cite{ChaKouMoy87}.  Orrick and Solomon showed that
appending a row and column to a specially normalized $20\times20$ Hadamard matrix on which a certain rank-1
update has been performed produces a
D-optimal design~\cite{OrrSol03}.  This construction is known as the {\em 3-normalized maximal
excess construction.}   The question presents itself: Do there exist optimal matrices that are
not derived from this construction?

To answer this, we observe that the condition $RR^{\mathrm{T}}=R^{\mathrm{T}}R=M$ imposes more structure on
$R$ than might be supposed.   Because of the structure of $M$, we segment $R$ into
a $(1+4+16)\times(1+4+16)$ matrix.  The row and column sums of the blocks formed
by this segmentation become relevant variables.
The conditions $M^2=RR^{\mathrm{T}}RR^{\mathrm{T}}=RMR^{\mathrm{T}}$
and
$RM=RR^{\mathrm{T}}R=MR$ imply a set of diophantine equations for these variables,
which are also subject to the obvious upper and lower bounds.  Solution of sets of equations of this type has been done in several places~\cite{MoyKou82,Orr05a}.  We omit the  proof of the following.
\begin{theo}
Let $r_j$ be the $j^{\text{th}}$ row of $R$ partitioned into $1+4+16$ elements as
$r_j=(a_j,b_j,c_j)$ and let $A_j$, $B_j$, and $C_j$ be the sums of the elements of $a_j$, $b_j$,
and $c_j$.  Then $(A_1,B_1,C_1)=(1,-4,16)$, $(A_j,B_j,C_j)=(-1,2,8)$ for $2\le j\le 5$, and
$(A_j,B_j,C_j)=(1,2,2)$ for $6\le j\le21$.  The columns of $R$ also have this structure.
\end{theo}

By the theorem, we can assume the first 5 rows take the form
\begin{equation*}
\mbox{
\tiny
\begin{math}
\begin{bmatrix}
+  & ---- & ++++ & ++++ & ++++ & ++++\\
\\
- &  -+++ & ---- & ++++ & ++++ & ++++\\
- &  +-++ & ++++ & ---- & ++++ & ++++\\
- &  ++-+ & ++++ & ++++ & ---- & ++++\\
- &  +++- & ++++ & ++++ & ++++ & ----\\
\end{bmatrix}
\end{math}
,
}
\end{equation*}
and the first 5 columns take this form transposed.  To determine the remaining $16\times16$ submatrix,
we partition it into a $4\times4$ array of $4\times4$ blocks.  The inner products of rows 6--21 with
rows 2--5, and the corresponding inner products for columns, force the row and column sums of the $4\times4$ blocks to take the values
\begin{equation*}
\begin{bmatrix}
2 & 0 & 0 & 0\\
0 & 2 & 0 & 0\\
0 & 0 & 2 & 0\\
0 & 0 & 0 & 2\\
\end{bmatrix}.
\end{equation*}

Without loss of generality we may now take the form of $R$ to be
\begin{equation}\label{form21}
\mbox{
\tiny
\begin{math}
\begin{bmatrix}
+  & ---- & ++++ & ++++ & ++++ & ++++\\
\\
- &  -+++ & ---- & ++++ & ++++ & ++++\\
- &  +-++ & ++++ & ---- & ++++ & ++++\\
- &  ++-+ & ++++ & ++++ & ---- & ++++\\
- &  +++- & ++++ & ++++ & ++++ & ----\\
\\
+ & -+++ & -+++\\
+ & -+++ & +-++\\
+ & -+++ & ++-+\\
+ & -+++ & +++-\\
\\
+ & +-++ & & -+++\\
+ & +-++ & & +-++\\
+ & +-++ & & ++-+\\
+ & +-++ & & +++-\\
\\
+ & ++-+ & & & -+++\\
+ & ++-+ & & & +-++\\
+ & ++-+ & & & ++-+\\
+ & ++-+ & & & +++-\\
\\
+ & +++- & & & & -+++\\
+ & +++- & & & & +-++\\
+ & +++- & & & & ++-+\\
+ & +++- & & & & +++-\\
\end{bmatrix}
\end{math}
.
}
\end{equation}
The missing blocks all have row and column sum 0.  Before proceeding to find them, we note that
if the first column is deleted, then row 1 together with rows 6--21 form an orthogonal set.  To complete
these 17 rows to a $20\times20$ Hadamard matrix, we need only add the 3 rows
\begin{equation*}
\mbox{
\tiny
\begin{math}
\begin{bmatrix}
--++ & ---- & ---- & ++++ & ++++\\
-+-+ & ---- & ++++ & ---- & ++++\\
-++- & ---- & ++++ & ++++ & ----\\
\end{bmatrix}
\end{math}
.
}
\end{equation*}
The implication is that any D-optimal matrix of order 21 is obtained from a $20\times20$ Hadamard
matrix by means of the 3-normalized maximal excess construction.  This construction amounts
to reversing the procedure by  replacing the 3 added rows with rows 2--5 of the above matrix minus
their first column, and finally appending column 1.  For details, see~\cite{OrrSol03}.  To get a
D-optimal matrix, the starting Hadamard matrix must have 3-normalized excess 76, which is the largest possible for order 20.  It turns out that any 3-normalized Hadamard
matrix of order 20  has row sums $(0_3,12,4_{16})$ and hence the maximal excess of 76.  We must
therefore consider all possible 3-normalizations of each of the 3 inequivalent $20\times20$ Hadamard matrices~\cite{Hal65}
in order to determine how many distinct D-optimal designs of order 21 there are.  The method we used
was the following:
\begin{itemize}
\item To each starting Hadamard matrix we applied the 3-normalized maximal excess construction.
\item We put the matrix in the form~(\ref{form21}) by suitable row and column permutations.
\item We further imposed the condition that the
$4\times 4$ block formed by the intersection of rows 6--9 and  columns 10--13 (the $(1,2)$ block) and 
the $4\times4$ block formed by the intersection of rows 14--17 and columns 18--21 (the $(3,4)$ block), both of which have row and column sums 0, be in one of two standard forms:
\begin{equation*}
\mbox{
\tiny
\begin{math}
\begin{bmatrix}--++\\-+-+\\+-+-\\++--\end{bmatrix}
\end{math}
,
}
\qquad
\mbox{
\tiny
\begin{math}
\begin{bmatrix}--++\\--++\\++--\\++--\end{bmatrix}
\end{math}
.
}
\end{equation*}
\item There remain two types of operation that preserve the overall structure we have imposed
on the matrix: (1) permutations of the
blocks that send diagonal blocks to diagonal blocks, coupled with the permutations of rows 2--5 and
columns 2--5 needed to preserve their form, and (2) automorphisms of 
the $(1,2)$ and $(3,4)$ blocks,  coupled with the transverse permutations needed to maintain
the form of the diagonal blocks.  We ran through these operations, stopping as soon as the resulting matrix  matched one on our canonical list (initially empty).  If all possible operations were exhausted
without finding a match, the current matrix was added to the canonical list.
\end{itemize}
By carrying out this procedure, we found 7 distinct forms.  These can be found at the website
{\em The Hadamard maximal determinant problem}~\cite{OrrSol05}.

\section{D-optimal designs of order 26}
There are two classes of D-optimal designs in order 26, relating to the two ways that 50 can be partitioned into two squares, $50=5^2+5^2$ and $50=7^2+1^2$.  Depending on which of
these we take, either all the matrices $A$, $B$, $C$, and $-D$ will have
row and column sums equal to 5 (the $(5,5)$-type), or $A$ and $-D$ will have row and column sums 7 and
$B$ and $C$ will have row and column sums 1 (the $(7,1)$-type).

Two standard constructions produce D-optimal designs in order 26.  Construction~1~\cite{Ehl64a} takes the submatrices
$A$ and $B$ to be circulant, and sets  $D=-A^{\mathrm{T}}$ and $C=B^{\mathrm{T}}$.  The
search for the initial rows of $A$ and $B$ can then be done exhaustively.  The complete set
of solutions, found by Yang~\cite{Yan68}, includes 3 matrices of $(5,5)$ type one of which is self-dual, the other of which form
a dual pair.  It also includes 3 matrices of $(7,1)$ type, again consisting of a self-dual matrix and
a dual pair.

Construction 2 (doubling) forms the matrices
\begin{equation*}
R=\begin{bmatrix}A & B\\A & -B\end{bmatrix}\qquad R'=\begin{bmatrix}A & A\\B & -B\end{bmatrix}
\end{equation*}
from maximal determinant matrices $A$ and $B$ of order 13~\cite{Woj64}.  $R$ and $R'$ are
necessarily of $(5,5)$ type.  Since there is a unique equivalence class in order 13, $A$ and
$B$ must be equivalent.  By performing certain transformations on $A$ or $B$, however, we
can obtain inequivalent matrices $R$.  Permuting or negating columns of $A$ or $B$ permutes 
or negates the corresponding columns of $R$
which does not change its equivalence class.  Negating row $j$ of $B$ interchanges
rows $j$ and $13+j$ of $R$.  Hence negating rows of $A$ or $B$ does not change the equivalence
class of $R$.  Permuting rows $j$ and $k$ of $B$ may, however change the equivalence class
of $R$.  We fix $A$, and consider the set of all possible permutations of rows of $B$.  By this
operation, and by transposition of $R$ we can produce exactly 367 equivalence classes.

Additional matrices have been produced in two ways.  One way is to use the gradient ascent
method of~\cite{OrrSolDowSmi03}. Good starting matrices for the ascent are $24\times 24$
Hadamard matrices, augmented with two random rows and columns, and maximal determinant
matrices of order 25, augmented with one random row and column.  Matrices of $(5,5)$ type
are produced with high probability, whereas $(7,1)$-type matrices are produced only a small
fraction of 1\% of the time.

The second way of producing new matrices is called {\em switching} and is extremely powerful.  Row switching
acts on 4 rows of the matrix.  These 4 rows must have the property that their Hadamard (element-wise)
product is the all 1 vector or its negation.  If this is the case, then the columns of the $4\times26$
matrix formed by extracting these rows will be of at most 8 different types.  Moreover, we consider a column and its negation to be the same type, so there are really only 4 different types.  Switching consists
of negating all columns of a given type.  It preserves D-optimality, but generally produces a matrix
inequivalent to the original matrix.  Negating any of the 4 different column types produces an equivalent result.  This method is equivalent to the method Denniston devised to enumerate $(25,9,3)$ designs~\cite{Den82}.  It was also applied to Hadamard matrices of order $8k$ in~\cite{Orr05b} which
contains proofs of many of the above statements.

As in~\cite{Orr05b} we define two matrices to be {\em Q-equivalent} if one can be obtained
from the other by some sequence of row permutation, row negation, and row switching,
combined with the corresponding column operations.  The associated equivalence classes 
are called {\em Q-classes.}  Since Q-equivalence is weaker than Hadamard equivalence,
we can regard the Q-class as composed of Hadamard equivalence classes.   A Q-class is
{\em self-dual} if it is identical to the set of the duals of the matrices contained in it.
By means of switching,
combined with equivalence checking using {\em nauty}~\cite{McK04,McK79} we have constructed the Q-classes of all $26\times26$ matrices obtained either by Constructions 1 and 2, or by gradient
ascent.

For $(5,5)$-type matrices, the dual pair from Construction 1 and all the matrices from Construction 2
are in the same self-dual Q-class, which has 8,545 elements.  That all matrices from Construction 2 (doubling) lie in
the same Q-class can be shown by the same method used to show the analogous property for
Hadamard matrices of order 32~\cite{Orr05b}.  The self-dual matrix from Construction 1 forms
a Q-class by itself.  There are 8 additional Q-classes which were constructed from matrices found by gradient ascent.  Three are singleton Q-classes, one self-dual,  two forming a dual pair.  Three are of size 4.  One is self-dual and the other two form a dual pair.  Two are
of size 5 and form a dual pair.  These 10 Q-classes together contain 8,571 Hadamard equivalence
classes.

For the $(7,1)$-type matrices, the three matrices from Construction 1 form singleton Q-classes, one
of which is self-dual, the other two forming a dual pair.  In addition, there is a self-dual Q-class of
size 1,310, constructed from matrices found by gradient ascent.  Combining these 1,313 matrices
with the 8,571 matrices of $(5,5)$-type, we find that there are at least 9,884 inequivalent D-optimal
designs of order 26.

The number of designs found in order 26 is surprising, considering that it is at least one or two orders
of magnitude higher than the numbers of designs in orders 24 (60 designs), 25 (78 designs), or 28 (487 designs).  It is tempting to speculate that this is accounted for by the existence of a doubling construction
for order 26.  Indeed, for Hadamard matrices, there appears to be a big jump in the number of
matrices at the orders in which doubling is possible, that is, the orders that are multiples of 8.  (Compare 3 matrices in order 20 with 60 matrices
in order 24, or 487 matrices in order 28 with more than 3,578,006 matrices in order 32.)  The matrices
produced by doubling are only a small fraction of the total, but they do provide abundant raw material
on which switching can act.

It is interesting in this connection to note that the D-optimal design in order 13 is obtained from the
Hadamard matrix of order 12 by the 3-normalized maximal excess construction.  (See
also~\cite{ConElkMar05}.)  One can formulate a doubled version of this construction which
will produce a D-optimal design of order 26 from a Hadamard matrix of order 24.  There is
considerable freedom when applying the construction which accounts for most of the combinatorial
explosion in the number of matrices generated.  All of the matrices produced in this way are
related to each other by switching but the converse is not true.  There are a small number of matrices
obtained by switching which cannot also be obtained by this doubled 3-normalized maximal
excess construction.

The considerations of the preceding paragraphs account for the large number of designs of $(5,5)$-type,
but the number of designs of $(7,1)$-type is also anomalously large.  We have no explanation for
this fact.  It is possible that there is an additional type of switching operation waiting to be found that converts matrices of one
type into the other type.

\section{Other applications of switching}
Switching is a very general technique.  
As noted above,
switching can also be applied to Hadamard matrices of order $n\equiv0\bmod 8$, and was originally
applied by Denniston in order 25.  A variant switching operation has been defined for Hadamard
matrices of order $n\equiv4\bmod8$~\cite{Orr05b}.

For orders $n\equiv2\bmod4$ we have, in addition to $n=26$,
applied it to $n=18$ where
we find that all three of Cohn's designs are related by switching.   We have also used it in orders
42 and 50 to produce vast numbers of matrices.  Note that there is no doubling construction in
orders 18 and 42.  On the other hand, one can easily prove the following negative result
\begin{propo}
Let $n\equiv6\bmod8$ and let $R$ be equivalent to a matrix whose Gram matrices satisfy
$M=M'=(n-2)I_n+2I_2\otimes J_{n/2}$ (the best possible form).  Then $R$ has no quadruple
of rows whose Hadamard product is the all-1 vector.
\end{propo}
This means that switching cannot be applied to matrices of the optimal form in orders congruent
to $6\bmod8$.  Nevertheless, we expect switching to apply to all orders $n\equiv2\bmod8$ for
which the best possible Gram matrix can be attained.  We have even applied switching in orders
such as 34 where the largest known determinant is suboptimal, and produced hundreds of
thousands of
inequivalent matrices.

Switching has proved to be useful in classification.
We have used switching in order 17 to show that all three D-optimal
designs are Q-equivalent and
in order 19 to show that Cohn's matrix and
the matrix found by Orrick and Solomon are Q-equivalent.  We also applied it to Tamura's record
determinant matrix in order 27 to produce 66 inequivalent matrices.  Curiously, we have never
found a D-optimal matrix whose order is congruent to $5$, $6$, or $7\bmod8$ to which switching
can be applied.

\section*{Acknowledgments}
I thank Bruce Solomon for many helpful conversations.
I used \emph{Mathematica} extensively during this
project and also thank Indiana University for the use of its  IBM RS/6000 SP computing platform.

\bibliographystyle{plain}
\bibliography{equivalenceRevised}

\end{document}